\documentclass[11pt]{article}
\usepackage{amsmath,amssymb,amsthm}

\usepackage{hyperref}
\hypersetup{
pdfauthor = {Boris Bukh},
pdftitle = {Non-trivial solutions to a linear equation in integers},
pdfsubject = {Mathematics},
pdfkeywords = {linear equation, trivial solutions, Sidon sets, Balog-Szemeredi-Gowers theorem},
pdfstartview = {FitH},
pdfpagemode = {None}}

\newtheorem{theorem}{Theorem}
\newtheorem{lemma}[theorem]{Lemma}
\newcommand*{\abs}[1]{\lvert #1\rvert}
\newcommand*{\bigabs}[1]{\left\lvert #1\right\rvert}
\newcommand*{\norm}[1]{\lVert #1\rVert}
\newcommand*{\Z}{\mathbb{Z}}

\author{Boris Bukh}
\title{Non-trivial solutions to a linear equation in integers}
\date{}
\begin{document}
\maketitle
\begin{abstract}
For $k\geq 3$ let $A\subset [1,N]$ be a set not containing
a solution to $a_1x_1+\dotsb+a_kx_k=a_1x_{k+1}+\dotsb+
a_kx_{2k}$ in distinct integers. We prove that
there is an $\varepsilon>0$ depending on the coefficients 
of the equation such that every 
such $A$ has $O(N^{1/2-\varepsilon})$ elements.
This answers a question of I.~Ruzsa.
\end{abstract}

\section*{Introduction}
The problems of estimating the size of a largest subset
of $[1,N]$ not containing a solution to a given linear
equation arise very frequently in combinatorial
number theory. For example, the sets with
no non-trivial solutions to $x_1+x_2-2x_3=0$ and
$x_1+x_2=x_3+x_4$ are sets with no arithmetic
progressions of length three, and Sidon sets respectively.
For several of the more prominent equations, like the two equations
above, there are large bodies of results that deal with
the structure and the size of solution-free sets. The first systematic
study of general linear equations was undertaken
by Ruzsa \cite{cite:ruzsa_lineareq_I,cite:ruzsa_lineareq_II}.
To ascribe a precise meaning to the concept of a ``set with no
non-trivial solution'' he introduced two definitions of
a trivial solution. One of them is that a solution is
non-trivial if all the variables are assigned different values.

For a fixed linear equation denote by $R(N)$ the size 
of a largest set of integers in $[1,N]$ with no solution
to the equation in distinct integers.
Ruzsa \cite{cite:ruzsa_lineareq_I} showed that if $k\geq 2$,
then for
the symmetric equation
\begin{equation}\label{theeq}
a_1x_1+\dotsc+a_kx_k=a_1x_{k+1}+\dotsb+a_kx_{2k}
\end{equation}
one has $R(N)=O(N^{1/2})$.
For
$k=2$, the estimate is tight for the Sidon equation.
Ruzsa gave the example of equation
\begin{equation}\label{ruzsaexample}
x_1+d(x_2+\dotsc+x_k)=x_{k+1}+d(x_{k+2}+\dotsc+x_{2k})
\end{equation}
and the set $A$ of integers 
 in $[1,N]$ whose development in the base
$d^2 k$ consists only of digits $0,\dotsc,d-1$.
Since addition of elements of $A$ in \eqref{ruzsaexample}
involves no carries in base $d^2 k$,
for every solution to \eqref{ruzsaexample} with
elements of $A$ one necessarily has that $x_1$ and $x_{k+1}$
are equal digit by digit. Thus $x_1=x_{k+1}$ implying
that $A$ contains
no solution \eqref{ruzsaexample} in distinct integers.
Since $\abs{A}=\Omega(N^{1/2-c/\log d})$,
this example shows that there is
no $\varepsilon>0$ such that the estimate 
$R(N)=O(N^{1/2-\varepsilon})$
holds for all equations of the form~\eqref{theeq}.
However, Ruzsa asked whether for $k\geq 3$ 
there is an $\varepsilon>0$ depending
on the coefficients of \eqref{theeq} such that 
$R(N)=O(N^{1/2-\varepsilon})$.
We answer in affirmative.
\begin{theorem}\label{thethm}
Let $\norm{a}_1=\sum_{i=1}^k \abs{a_i}$. If
$k\geq 3$ and $a_i\neq 0$ for $1\leq i\leq k$,
then 
\begin{equation*}
R(N)=O\bigl(N^{\frac{1}{2}-\frac{1}{c(k)\norm{a}_1}}\bigr).
\end{equation*}
\end{theorem}

To explain the idea behind the proof of theorem~\ref{thethm}
we examine the proof that $R(N)=O(N^{1/2})$. Let
$a_{k+i}=-a_i$. Then the equation \eqref{theeq} can
be written as $\sum a_i x_i=0$.
Assume that $A\subset [1,N]$ has $M$ elements and 
contains only trivial solutions to~\eqref{theeq}.
Let $r(m)$ denote the number of solutions to
\begin{equation*}
a_1 x_1+\dotsb+a_k x_k=m,\qquad x_i\in A.
\end{equation*}
Then $E=\sum r(m)^2$ is the total number of solutions
to~\eqref{theeq} in~$A$. Since
there are at most $\norm{a}_1 N$ 
values of $m$ for which $r(m)>0$ 
by Cauchy-Schwarz it follows that 
\begin{equation}\label{trivlowerbound}
E=\sum r(m)^2\geq \frac{M^{2k}}{\norm{a}_1 N}.
\end{equation}
The next step is to bound the number of solutions with $x_i=x_j$
from above. Let $1\leq r\leq 2k$ be any index other than $i$
or~$j$. Since 
\begin{equation}\label{linearsum}
x_r=-\frac{1}{a_r}\Bigl(\sum_{l\neq i,j,r} a_l x_l+(a_i+a_j)x_i\Bigr),
\end{equation}
it follows that equation \eqref{theeq} uniquely determines $x_r$
in terms of other $x$'s.  Thus at most one in every
$M$ assignments of variables results in a solution to~\eqref{theeq}.
Since there $\binom{2k}{2}$ choices of $i$ and $j$
\begin{equation}\label{trivupperbound}
E\leq \binom{2k}{2}\frac{1}{M}M^{2k-1}.
\end{equation}
The bound $\abs{A}=O(N^{1/2})$ follows by comparison  
with~\eqref{trivlowerbound}.

In order for this argument to be tight the linear form
on the right side of \eqref{linearsum} should often take
values in $A$. Heuristically it corresponds to the sumset
$\sum_{l\neq i,j,r} x_l\cdot A+(x_i+x_j)\cdot A$ being not much larger than
$A$ itself, where $t\cdot A=\{ta : a\in A\}$ is the $t$-dilate
of~$A$. Pl\"unnecke-Ruzsa inequalities tells us that
if some sumset involving $A$ is not much larger than $A$,
then every sumset involving only $A$ is not much larger 
than~$A$. In particular we expect $\sum_{l\leq k} x_1\cdot A$
to be of size only $M$, as opposed to $\sim N$ which 
was assumed in the derivation of the lower bound \eqref{trivlowerbound}.
Therefore, we expect that if \eqref{trivupperbound} is
close to being tight, then \eqref{trivlowerbound} is not,
and vice versa.
The remainder of the paper is a rigorous justification
of the heuristic argument above.

\section*{Proof}
This section is organized as follows. First we state the tools from 
additive combinatorics that we need. Then we introduce the notation
that is going to be used in the proof of theorem~\ref{thethm}. 
Finally, after a couple of preliminary lemmas the theorem~\ref{thethm}
is proved.

\begin{lemma}[Ruzsa's triangle inequality \cite{cite:tao_vu_book}, lemma~2.6]
For any finite $A,B,C\subset \Z$ we have 
\begin{equation*}
\abs{A-C}\leq\frac{\abs{A-B}\abs{B-C}}{\abs{B}}.
\end{equation*}
\end{lemma}
\begin{lemma}[Pl\"unnecke's inequality \cite{cite:tao_vu_book}, corollary~6.26]
For any finite sets $A,B\subset \Z$ if $\abs{A+B}\leq K\abs{A}$, then
$\abs{kB}\leq K^k\abs{A}$.
\end{lemma}
We will also use the hypergraph version of Balog-Szemer\'edi-Gowers theorem
due to Sudakov, Szemer\'edi and Vu \cite{cite:hypergraph_balog_szem}. 
Let $A_1,\dotsc,A_k$ be sets of integers. If $H$ is a subset of
$A_1\times\dotsb A_k$, then $\sum_H A_i$ is the collection of all
sums $a_1+\dotsb+a_k$, where $(a_1,\dotsc,a_k)\in H$.
\begin{lemma}[\cite{cite:hypergraph_balog_szem}, theorem~4.3]\label{lemmabalogszem}
For any integer $k\geq 1$ there are positive-valued polynomials 
$f_k(x,y)$ and $g_k(x,y)$ such that the following holds. Let $n,C,K$
be positive numbers. If $A_1,\dotsc,A_k$ are sets of $n$ positive integers,
$H\subset A_1\times \dotsb\times A_k$ with $\abs{H}\geq n^k/K$ and 
$\bigabs{\sum_H A_i}\leq Cn$, then one can find subsets $A_i'\subset A_i$
such that
\begin{align*}
&\abs{A_i'}\geq \frac{n}{f_k(C,K)},\qquad\text{for all }1\leq i\leq k,\\
&\abs{A_1'+\dotsb+A_k'}\leq g_k(C,K)n.
\end{align*}
\end{lemma}

From the heuristic argument above it is clear
that we will need to count number of solutions
of various equations. It is therefore advantageous
to introduce appropriate notation. Let
$r(A_1,\dotsc,A_k;m)$ denote the number of
solutions to
\begin{equation*}
a_1+\dotsb+a_k=m,\qquad a_i\in A_i.
\end{equation*}
Then $E(A_1,\dots,A_k;B_1,\dotsc,B_l)=\sum_m r(A_1,\dotsc,A_k;m)r(B_1,\dotsc,
B_l;m)$ counts the number of solutions to
\begin{equation*}
a_1+\dotsb+a_k=b_1+\dotsb+b_l,\qquad a_i\in A_i,\ b_j\in B_j.
\end{equation*}
We write $E(A_1,\dotsc,A_k)$ to denote $E(A_1,\dotsc,A_k;A_1,\dotsc,A_k)$.
\begin{lemma}
There is a positive-valued polynomial $h_k(x)$ such that
if $A_1,\dotsc,A_k$ are $k$ sets of integers with $n$ elements each,
and $E(A_1,\dotsc,A_k)\geq c n^{2k-1}$, then there are subsets $A_i'\subset A_i$
such that
\begin{align*}
&\abs{A_i'}\geq h_k(c)n,\qquad\text{for all }1\leq i\leq k,\\
&\abs{A_1'+\dotsb+A_k'}\leq \frac{n}{h_k(c)}.
\end{align*}
\end{lemma}
\begin{proof}
For brevity of notation write $E=E(A_1,\dotsc,A_k)$
and $r(m)=r(A_1,\dotsc,A_k;m)$. Let $S=\{m : r(m)>c n^{k-1}/2\}$, and
$H=\{(a_1,\dotsc,a_k)\in A_1\times\dotsb\times A_k :
a_1+\dotsb+a_k\in S\}$. Note that $\sum_H A_i=S$ and 
\begin{equation*}
\abs{S} \leq \frac{n^k}{cn^{k-1}/2}=\frac{2}{c}n.
\end{equation*}
Since $r(m)\leq n^{k-1}$ for every $m$, 
and $\sum_{m\not\in S}r(m)=n^k-\abs{H}$,
it follows that
\begin{align*}
E=\sum_{m\in S} r(m)^2+\sum_{m\not\in S} r(m)^2
\leq \abs{H}n^{k-1}+(n^k-\abs{H})cn^{k-1}/2.
\end{align*}
Therefore $\abs{H}\geq cn^k/2$, and we deduce
from lemma \ref{lemmabalogszem} the existence of the subsets
$A_i'$ with the desired properties.
\end{proof}
\begin{lemma}\label{lemdilate}
For every $k\geq 2$ there is a positive-valued polynomial $p_k(x)$ such
if $t_1,\dotsc,t_k$ are any $k$ positive integers, and $A\subset \Z$ is
an $n$-element set satisfying $E(t_1\cdot A,\dotsc,t_k\cdot A)\geq c n^{2k-1}$,
then for every $l$-tuple of positive integers $s_1,\dotsc,s_l$ we have
$E(s_1\cdot A,\dotsc,s_l\cdot A)\geq p_k(c)^{\norm{s}_1} n^{2l-1}$.
\end{lemma}
\begin{proof}
Apply the lemma above with $A_i=t_i\cdot A$ to obtain $A_i'$.
Ruzsa's triangle inequality applied to
$-A_1'$ and $A_2'+\dotsb+A_k'$ yields that 
\begin{equation*}
\abs{A_1'-A_1'}\leq
\frac{\abs{A_1'+A_2'+\dotsb+A_k'}^2}{\abs{A_2'+\dotsb+A_k'}}
\leq \frac{n}{h_k(c)^3}.
\end{equation*}
Set $\tilde{A}=(1/t_1)\cdot A_1'$. Then $\tilde{A}$ is
a subset of $A$ satisfying $\abs{\tilde{A}}\geq
h_k(c)n$ and $\abs{\tilde{A}-\tilde{A}}\leq n/h_k(c)^3$.
Pl\"unnecke's inequality then implies that 
$\abs{\norm{s}_1\tilde{A}}\leq h_k(c)^{-3\norm{s}_1}n$. 
The inclusion $s_1\cdot \tilde{A}+\dotsb+s_l\cdot \tilde{A}
\subset \norm{s}_1\tilde{A}$ together with Cauchy-Schwarz
inequality yield
\begin{align*}
E(s_1\cdot A,\dotsc,s_l\cdot A)&\geq E(s_1\cdot \tilde{A},\dotsc,s_l\cdot \tilde{A})\\
&\geq \frac{\abs{\tilde{A}}^{2l}}{\abs{s_1\cdot \tilde{A}+\dotsb+s_l\cdot \tilde{A}}}\\
&\geq \frac{h_k(c)^{2l}n^{2l}}{\abs{\norm{s}_1 \tilde{A}}}\\
&\geq h_k(c)^{3\norm{s}_1+2l} n^{2l-1}.
\end{align*}
Since $s_1,\dotsc,s_l$ are positive integers, $\norm{s}_1\geq l$,
and the lemma follows.
\end{proof}
\begin{proof}[Proof of theorem~\ref{thethm}]
Without loss of generality we may assume that
$a_1,\dotsc,a_k$ are positive.
Rewrite equation \eqref{theeq} as
\begin{equation}\label{rewritetheeq}
a_1x_1+\dotsc+a_{2k}x_{2k}=0,
\end{equation}
where $a_{k+i}=-a_i$. Let $A\subset [1,N]$ be an $M$-element set 
containing 
only trivial solutions to~\eqref{rewritetheeq}.
The number of solutions to equation~\eqref{rewritetheeq}
with variables in $A$ is $E=E(a_1\cdot A,\dotsc,a_k\cdot A)$. Let
$T_{i,j}$ be the number of solutions with $x_i=x_j$. By the
pigeonhole principle for at least one pair $i\neq j$
we have $T_{i,j}\geq E/\binom{2k}{2}$. Fix that pair.

Next we partition the index set $\{1,\dotsc,2k\}\setminus
\{i,j\}$ arbitrarily into sets $I_1=\{l_{1,1},\dotsc,l_{1,k-1}\}$ and 
$I_2=\{l_{2,1},\dotsc,l_{2,k-1}\}$ with $k-1$ 
elements each. For convenience write 
\begin{equation*}
r_1(m)=r(a_{l_{1,1}}\cdot A_{l_{1,1}},\dotsc,a_{l_{1,k-1}}\cdot A_{l_{1,k-1}};m)
\end{equation*} 
and $E_1=\sum r_1(m)^2$. Define $r_2$ and $E_2$ analogously with respect 
to~$I_2$.
Then by Cauchy-Schwarz
\begin{equation*}
T_{i,j}=\sum_{x\in A} \sum_m r_1(m)r_2\bigl(-(a_i+a_j)x-m\bigr)
\leq \sum_{x\in A} E_1^{1/2}E_2^{1/2}=M E_1^{1/2}E_2^{1/2}.
\end{equation*}
Either $E_1$ or $E_2$ is at least $E/\binom{2k}{2}M$. We can assume it
is true of $E_1$. Then by lemma~\ref{lemdilate} above we have
\begin{equation*}
E\geq p_{k-1}\left(\frac{E_1}{M^{2k-3}}\right)^{\norm{a}_1}M^{2k-1}
 \geq p_{k-1}\left(\frac{E}{\binom{2k}{2}M\, M^{2k-3}}\right)^{\norm{a}_1} M^{2k-1}.
\end{equation*}
If $\deg p_{k-1}=d$, then we obtain that
\begin{equation*}
E=O(M^{2k-2-1/(d\norm{a}_1-1)}).
\end{equation*}
Since $E\geq M^{2k}/\norm{a}_1 N$, the theorem follows.
\end{proof}

\section*{Conclusion}
Theorem~\ref{thethm} gives the estimate
$R(N)=O\bigl(N^{\frac{1}{2}-\frac{1}{r}}\bigr)$ with
$r=c(k)\norm{a}_1$.
Ruzsa's example shows that no estimate
better than $r=c(k)\log \norm{a}_1$
can be true. 

After this paper was written, a Pl\"unnecke-type
estimate on sums of dilates of the form
$s_1\cdot A+\dotsb+s_l\cdot A$, which appears
in the proof of lemma~\ref{lemdilate}, was established
in \cite{cite:bukh_sumsodilates}. Whereas Pl\"unnecke's inequality
and the inclusion $s_1\cdot A+\dotsb+s_l\cdot A\subset
\norm{s}_1 A$ yield the exponent of $\norm{s}_1$, the new
inequality yields the exponent $C\log \norm{s}_1$ for
an absolute constant~$C$. Plugging that into the proof
of lemma~\ref{lemdilate}, one obtains that
theorem~\ref{thethm} is valid with $r=c(k)\log\norm{a}_1$,
which is sharp in view of Ruzsa's example.

\vspace{0.3cm}\noindent{\bf  Acknowledgement.}\, I thank Benjamin Sudakov for stimulating discussions.

\bibliographystyle{alpha}
\bibliography{ruzsatriv}

\begin{thebibliography}{SSV05}

\bibitem[Buk07]{cite:bukh_sumsodilates}
Boris Bukh.
\newblock Sums of dilates.
\newblock \hhref{0711.1610}, Nov 2007.

\bibitem[Ruz93]{cite:ruzsa_lineareq_I}
Imre~Z. Ruzsa.
\newblock Solving a linear equation in a set of integers. {I}.
\newblock {\em Acta Arith.}, 65(3):259--282, 1993.

\bibitem[Ruz95]{cite:ruzsa_lineareq_II}
Imre~Z. Ruzsa.
\newblock Solving a linear equation in a set of integers. {II}.
\newblock {\em Acta Arith.}, 72(4):385--397, 1995.

\bibitem[SSV05]{cite:hypergraph_balog_szem}
B.~Sudakov, E.~Szemer{\'e}di, and V.~H. Vu.
\newblock On a question of {E}rd{\H o}s and {M}oser.
\newblock {\em Duke Math. J.}, 129(1):129--155, 2005.
\newblock \url{http://www.math.princeton.edu/~bsudakov/erdos-moser.pdf}.

\bibitem[TV06]{cite:tao_vu_book}
Terence Tao and Van~H. Vu.
\newblock {\em Additive Combinatorics}, volume 105 of {\em Cambridge studies in
  advanced mathematics}.
\newblock CUP, 2006.

\end{thebibliography}

\end{document}